\documentclass[12pt, reqno]{amsart}
\usepackage[utf8]{inputenc}
\usepackage{amssymb,mathtools,cite,enumerate,color,eqnarray,hyperref,amsfonts,amsmath,amsthm,setspace,tikz,verbatim,times}

\usepackage[a4paper,top=2cm,bottom=2cm,left=2.2cm,right=2.2cm]{geometry}

\usepackage[T1]{fontenc}
\usepackage[greek,english]{babel}

\numberwithin{equation}{section}

\definecolor{ao(english)}{rgb}{0.0, 0.5, 0.0}

\hypersetup{colorlinks=true, linkcolor=ao(english),citecolor=ao(english)}

\usepackage[normalem]{ulem}

\usepackage{color, xcolor}

\newtheorem{theorem}{Theorem}[section]
\newtheorem{conjecture}{Conjecture}[section]
\newtheorem{definition}{Definition}[section]
\newtheorem{corollary}{Corollary}[section]
\newtheorem{lemma}{Lemma}[section]
\newtheorem{remark}{Remark}[section]

\title[Arithmetic Properties of $k$-tuple $\ell$-regular Partitions]{Arithmetic Properties of $k$-tuple $\ell$-regular Partitions}

\author[H. Nath]{Hemjyoti Nath}
\address[H. Nath]{Lokhra chariali, Guwahati 781040, Assam, India}
\email{hemjyotinath40@gmail.com}

\author[M. P. Saikia]{Manjil P. Saikia}
\address[M. P. Saikia]{Mathematical and Physical Sciences division, School of Arts and Sciences, Ahmedabad University, Ahmedabad 380009, Gujarat, India}
\email{manjil@saikia.in}

\author[A. Sarma]{Abhishek Sarma}
\address[A. Sarma]{Department  of Mathematical Sciences, Tezpur University, Napaam,  Tezpur 784028, Assam, India}
\email{abhitezu002@gmail.com}



\date{}

\begin{document}

\begin{abstract}
In this paper, we study arithmetic properties satisfied by the $k$-tuple $\ell$-regular partitions. A $k$-tuple of partitions $(\xi_1, \xi_2, \ldots, \xi_k)$ is said to be $\ell$-regular if all the $\xi_i$'s are $\ell$-regular. We study the cases $(\ell, k)=(2,3), (4,3), (\ell, p)$, where $p$ is a prime, and even the general case when both $\ell$ and $k$ are unrestricted. Using elementary means as well as the theory of modular forms we prove several infinite family of congruences and density results for these family of partitions.
\end{abstract}

\maketitle

\vspace{5mm}

\section{Introduction}

A partition $\lambda$ of a natural number $n$ is a nonincreasing sequence of natural numbers whose sum is $n$. If $\lambda=(\lambda_1, \lambda_2, \ldots, \lambda_k)$ such that $\lambda_1\geq \lambda_2\geq \cdots \geq \lambda_k$ and $\sum\limits_{i=1}^k \lambda_i=n$, then $\lambda$ is called a partition of $n$, and $\lambda_i$'s are called the parts of the partition $\lambda$. For instance, the $7$ partitions of $5$ are
\[
5, 4+1, 3+2, 3+1+1, 2+2+1, 2+1+1+1, 1+1+1+1+1.
\]
We denote by $p(n)$ the number of partitions of $n$, and its generating function was given by Euler to be
\[
\sum_{n\geq 0}p(n)q^n=\frac{1}{\prod_{i=1}^\infty(1-q^i)}.
\]
For ease of notation, we write $(a;q)_\infty:=\prod\limits_{i=0}^\infty(1-aq^i)$ and $f_k:=(q^k;q^k)_\infty$. Thus, Euler's generating function becomes
\[
\sum_{n\geq 0}p(n)q^n=\frac{1}{(q;q)_\infty}=\frac{1}{f_1}.
\]

Partitions have been studied since the time of Euler, and several well-known mathematicians have explored their properties. Prominent among them is Ramanujan, who in 1920 \cite{Ramanujan} proved the following amazing congruences that the partition function satisfies: for all $n\geq 0$, we have
\begin{align*}
    p(5n+4)&\equiv 0\pmod 5,\\
 p(7n+5)&\equiv 0\pmod 7,\\
     p(11n+6)&\equiv 0\pmod{11}.
\end{align*}
Since then, one strand of research related to partitions is to find such Ramanujan-type congruences for partitions as well as for generalized partitions.
For a general overview of the area of partitions, we refer the reader to the excellent books by Andrews \cite{gea1} and Johnson \cite{john}.

Among the class of generalized partitions, a frequently studied class is that of $\ell$-regular partitions, for $\ell>1$. By an $\ell$-regular partition of $n$ we mean a partition of $n$ where no parts are divisible by $\ell$. Let $b_\ell(n)$ denote the number of $\ell$-regular partitions of $n$, then we have the following generating function
\[
\sum_{n\geq 0}b_\ell(n)q^n=\frac{f_\ell}{f_1}.
\]

In this paper, we are interested in a more general class of partitions, which we call $k$-tuple $\ell$-regular. A partition $k$-tuple $(\xi_1, \xi_2, \ldots, \xi_k)$ is called a $k$-tuple $\ell$-regular partition if all of the $\xi_i$'s are themselves $\ell$-regular partitions. Let us denote the number of such partitions of $n$ by $T_{\ell,k}(n)$. It is easy to see that its generating function is given by
\begin{equation}\label{eq:gf-lk}
    \sum_{n\geq 0}T_{\ell,k}(n)q^n=\dfrac{f_\ell^k}{f_1^k}.
\end{equation}
When $k=3$, we suppress the value of $k$ and just use the notation $T_{\ell,3}(n)=T_\ell(n)$. So, we get
\begin{equation}\label{e1.0.0.0}
    \sum_{n\geq 0} T_{\ell}(n)q^n = \dfrac{f_\ell^3}{f_1^3}.
\end{equation}

Although, $\ell$-regular partitions are very well studied, it seems that $k$-tuple $\ell$-regular partitions have not received the same attention. In this paper, we remedy this situation and study various arithmetic properties that the $T_{\ell, k}(n)$ function satisfies. The case when $\ell=k=3$ was first studied by Adiga and Dasappa \cite{AdigaDasappa}, the case $\ell=3$ and $k=9, 27$ were studied by Baruah and Das \cite{BaruahDas}, the case $\ell=3, k=6$ was studied by Murugan and Fathima \cite{MuruganFathima}, and very recently Nadji and Ahmia \cite{NadjiAhmia} studied the cases \(\ell=2, k=3\) and $\ell=k=3$. Here, we not only study the cases \(\ell=2, k=3\) and $\ell=k=3$, extending some of the results of Nadji and Ahmia \cite{NadjiAhmia}, but also the cases $(\ell, k)=(4,3), (\ell, p)$, for a prime $p$ as well as the more general case when $\ell$ and $k$ are unrestricted. Our proof techniques come from both elementary means as well as from the theory of modular forms.

We begin our results by first proving a general congruence that $T_{\ell,p}(n)$ satisfies, where $p$ is a prime. The proof is short and simple, so we complete it here. 
\begin{theorem}
    Let $p$ be a prime and $l$ be a non-negative integer. Then
    \begin{align}
        T_{\ell,p}(pn+r)\equiv 0 \pmod p\label{cong:0 mod p} 
    \end{align}
     for $r\in\{1,2,\ldots, p-1\}$.
\end{theorem}
\begin{proof}
    Putting $k = p$ in \eqref{eq:gf-lk}, we have
    \begin{align*}
        \sum_{n\geq 0}T_{\ell, p}(n)q^n&=\dfrac{f_\ell^p}{f_1^p}\equiv\dfrac{f_{\ell p}}{f_p}\pmod p.
    \end{align*}
Comparing the coefficients of $pn+r$ for $r\in\{1,2,\ldots, p-1\}$ on both sides, we arrive at \eqref{cong:0 mod p}.
\end{proof}
\noindent In the above proof, we have used the following easily verifiable identity: for a prime $p$, and positive integers $k$ and $l$, we have
     \begin{align}\label{e0.1}
         f_{k}^{p^l} \equiv f_{pk}^{p^{l-1}} \pmod{p^l}.
     \end{align}
We will use this fact without commentary in the sequel.

Before proceeding to our other results, we state the following result without proof, which follows very easily from an application  of \eqref{e2.0.3.3} and \eqref{e0.2}, stated in the next section.

\begin{theorem}\label{t0.1}
    For $n\geq0$, let $T_n$ be the $n$-th triangular number, then
    \begin{equation*}
        T_{2}(9n+1) = \begin{cases}
            3 \pmod{6} \hspace{1mm} \text{if} \quad n = T_n,\\
            0 \pmod{6} \hspace{1.5mm} \text{otherwise}.
        \end{cases}
    \end{equation*}
\end{theorem}

The next few results give several infinite family of congruences for $T_{\ell}(n)$ when $\ell=2,4$.
\begin{theorem}\label{c1.4}
    For all $n\geq 0$ and $\alpha\geq 0$, we have
 \begin{align}     
T_{2}\left(3^{4\alpha+2}n+\sum_{i=0}^{2\alpha}3^{2i}+3^{4\alpha+1}\right)&\equiv 0\pmod{24}, \label{c0.1.4}\\
T_{2}\left(3^{4\alpha+2}n+\sum_{i=0}^{2\alpha}3^{2i}+2\cdot 3^{4\alpha+1}\right)&\equiv 0\pmod{24}, \label{c1.1.4}\\
T_{2}\left(3^{4\alpha+4}n+\sum_{i=0}^{2\alpha+1}3^{2i}+3^{4\alpha+3}\right)&\equiv 0\pmod{24}, \label{c2.1.4}\\
T_{2}\left(3^{4\alpha+4}n+\sum_{i=0}^{2\alpha+1}3^{2i}+2\cdot 3^{4\alpha+3}\right)&\equiv 0\pmod{24}. \label{c3.1.4}
 \end{align}
\end{theorem}

\begin{remark}
    Nadji and Ahmia \cite[Theorem 3]{NadjiAhmia} proved the above congruences modulo $12$, here we generalize them to modulo $24$.
\end{remark}

\begin{theorem}\label{c1.4.1}
    For all $n\geq 0$ and $\alpha \geq 0$, we have
\begin{align}
        T_{4}\left( 3^{2\alpha +2 }n + \dfrac{17 \cdot 3^{2\alpha+1}-3}{8} \right) & \equiv 0 \pmod{3}, \label{e3.0}\\
        T_{4}\left( 3^{2\alpha +3 }n + \dfrac{19 \cdot 3^{2\alpha+2}-3}{8} \right) & \equiv 0 \pmod{3}, \label{e3.1}\\
        T_{4}\left( 27 \cdot 5^{2\alpha}n + \dfrac{171 \cdot 5^{2\alpha}-3}{8} \right) & \equiv 0 \pmod{3}. \label{e2.9}
\end{align}
\end{theorem}

Theorem \ref{c1.4} is proved in Section \ref{sec:pf-1} and Theorem \ref{c1.4.1} is proved in Section \ref{sec:pf-2}. A corollary of Theorem \ref{c1.4.1} is also given in Section \ref{sec:pf-2}.

For the next few results, we will need the Legendre symbol, which for a prime $p\geq3$ is defined as
\begin{align*}
\left(\dfrac{a}{p}\right)_L:=\begin{cases}\quad1,\quad \text{if $a$ is a quadratic residue modulo $p$ and $p\nmid a$,}\\\quad 0,\quad \text{if $p\mid a$,}\\~-1,\quad \text{if $a$ is a quadratic nonresidue modulo $p$.}
\end{cases}
\end{align*}

\begin{theorem}\label{t0.1.0.0}
    Let $p \geq 3$ be prime and let $r$, $1 \leq r \leq p-1$, be such that $\left( \dfrac{8r+1}{p} \right)_L = -1$, then for all $n \geq 0$, we have $T_2(9(pn+r)+1) \equiv 0 \pmod{6}$.
\end{theorem}

\begin{theorem}\label{t0.0.1}
    Let $p \geq 5$ be a prime. If $\left(\dfrac{-2}{p}\right)_L=-1$, then for $n,\alpha \geq 0$ with $p \nmid n$,
    \begin{equation}\label{e50.0}
        T_{2}\left( 9p^{2\alpha+1}n + \dfrac{9p^{2\alpha+2}-1}{8} \right) \equiv 0 \pmod{6}.
    \end{equation}
\end{theorem}

\begin{theorem}\label{t2}
    Let $a(n)$ be defined by
    \begin{equation*}
        \sum_{n=0}^{\infty}a(n)q^n = \dfrac{f_3^6}{f_1}.
    \end{equation*}
Let $p\geq5$ be a prime, we define
\begin{equation*}
    \omega(p) := a\left( \dfrac{17}{24}(p^2-1) \right)+p\left( \dfrac{2}{p} \right)_L\left( \dfrac{\frac{-17}{24}(p^2-1)}{p} \right)_L.
\end{equation*}
    \begin{enumerate}[(i)]
         \item If $\omega(p) \equiv 0 \pmod{2}$, then for $n,k\geq 0$ and $1\leq j \leq p-1$, we have
    \begin{equation}\label{t3.1}
            T_{2}\left( 3p^{4k+4}n + 3p^{4k+3}j + \dfrac{17\cdot p^{4k+4}-1}{8} \right) \equiv 0 \pmod{12}.
        \end{equation}
Furthermore, if $17 \nmid (24n + 17)$ and $n \not\equiv 0 \pmod{17}$, then for all $k \geq 0$, we have
    \begin{equation}\label{t3.3.1}
             T_{2}(3n+2) \equiv T_{2}\left( 3\cdot 17^{4k+2}n + \dfrac{17^{4k+3}-1}{8} \right) \pmod{12}.
        \end{equation}

        \item If $\omega(p) \equiv 1 \pmod{2}$, then for $n,k\geq 0$ and $1\leq j \leq p-1$, we have
      \begin{equation}\label{t2.2}
         T_{2}\left( 3p^{6k+6}n + 3p^{6k+5}j + \dfrac{17\cdot p^{6k+6}-1}{8} \right) \equiv 0 \pmod{12}.
         \end{equation}
Furthermore, if $p \nmid (24n + 17)$ and $n \not\equiv -17 \cdot 24^{-1} \pmod{p}$, then for all $k \geq 0$, we have
    \begin{equation}\label{t3.3}
             T_{2}\left( 3p^{6k+2}n + \dfrac{17\cdot p^{6k+2}-1}{8} \right) \equiv 0 \pmod{12}.
        \end{equation}
\end{enumerate}      
\end{theorem}

\begin{remark}
As an application of Theorem \ref{t2} we mention a congruence that we can derive: Using \textit{Mathematica}, we find that $a(17) \equiv 1 \pmod{2}$. Using Theorem \ref{t2} with $p=5$, we see that $\omega(5) \equiv 0 \pmod{2}$ and hence by \eqref{t3.1}, we obtain that for any integer $n, k \geq 0$ and $j=1$, we have \[T_2\left(1875\cdot 5^{4k}n + 375\cdot 5^{4k}+\dfrac{10625\cdot 5^{4k}-1}{8}\right) \equiv 0 \pmod{12}.\] Thus, Theorem \ref{t2} gives us infinite family of congruences.
\end{remark}

We recall the Dedekind $\eta$- function, given by
\begin{equation}\label{eta 1}
    \eta(z) := q^{1/24}f_1,
\end{equation}
where $q=e^{2\pi iz}$ and $z$ lies in the complex upper half plane $\mathbb{H}$. The well known $\Delta$-function is denoted by
\begin{equation*}
    \Delta(z) := \eta(z)^{24} = \sum_{n=1}^{\infty}\tau(n)q^n.
\end{equation*}

\begin{theorem}\label{t4}
Let $p$ be an odd prime. Then the following statements hold:
\begin{enumerate}[(i)]
    \item Suppose that $s$ is an integer satisfying $1\leq s \leq 8p$, $s \equiv 1 \pmod{8}$ and $\left( \dfrac{s}{p} \right)_L = -1$. Then,
\begin{equation}\label{e12.0.1}
    T_{2}\left( pn+\dfrac{s-1}{8} \right) \equiv 0 \pmod{2}.
\end{equation}
If $\tau(p) \equiv 0 \pmod{2}$, then, for all $n \geq 0, k\geq 1$,
\begin{equation}\label{e12.0.2}
    T_{2}\left( p^{2k+1}n+\dfrac{sp^{2k}-1}{8} \right) \equiv 0 \pmod{2}.
\end{equation}
    \item Suppose that $r$ is an integer  such that $1\leq r \leq 8p$, $rp \equiv 1 \pmod{8}$ and $(r,p)=1$. If $\tau(p) \equiv 0 \pmod{2}$, then, for all $n \geq 0, k\geq 1$,
    \begin{equation}\label{e12.0.3}
    T_{2}\left( p^{2k+2}n + \dfrac{rp^{2k+1}-1}{8} \right) \equiv 0 \pmod{2}.
\end{equation}
\end{enumerate}
\end{theorem}

We prove Theorem \ref{t0.1.0.0} in Section \ref{sec:thmn2}, Theorem \ref{t0.0.1} in Section \ref{s4}, Theorem \ref{t2} in Section \ref{s5}, and Theorem \ref{t4} in Section \ref{s3}.

Using the theory of modular forms, we prove the following result in Section \ref{sec:thmn1}.

\begin{theorem}\label{thm1.00}
	Let $k, n$ be nonnegative integers. For each $i$ with $1\leq i \leq k+1$, if $p_i \geq 5$ is a prime such that $p_i \not\equiv 1 \pmod 8$, then for any integer $j \not\equiv 0 \pmod {p_{k+1}}$
	\begin{align*} 
	T_2\left(9p_1^2\dots p_{k+1}^2n + \frac{9p_1^2\dots p_{k}^2p_{k+1}(8j+p_{k+1})-1}{8}\right) \equiv 0 \pmod{6}.
\end{align*}
\end{theorem}
As a consequence of the above theorem, the following corollary can now be easily deduced.
\begin{corollary}
    Let $p\geq 5$ be a prime such that $p \not\equiv 1\pmod{8}$. By taking all the primes $p_1, p_2, \ldots, p_{k+1}$ to be equal to the same prime $p$ in Theorem \ref{thm1.00}, 
we obtain the following infinite family of congruences
	\begin{align*} 
	T_2\left( 9p^{2(k+1)}n + 9p^{2k+1}j + \frac{9 p^{2(k+1)}-1}{8}\right) \equiv 0 \pmod 6,
	\end{align*}
	where  $j \not\equiv 0 \pmod p$.
\end{corollary} 

As an example,  taking $p=5$, we have that for all $n, k\geq 0$ and $j\not\equiv 0\pmod{5}$, we have
	\begin{align*} 
	T_2\left(225\cdot 5^{2k}n + 45\cdot 5^{2k}j + \frac{225\cdot 5^{2k}-1}{8}\right) \equiv 0 \pmod{6}.
	\end{align*}

Along with the study of Ramanujan-type congruences, the study of distribution of the coefficients of a formal power series modulo $M$ is also an interesting problem to explore. Given an integral power series $A(q) := \displaystyle\sum_{n=0}^{\infty}a(n)q^n $ and $0 \leq r \leq M$, the arithmetic density $\delta_r(A,M;X)$ is defined as
\begin{equation*}
       \delta_r(A,M;X) = \frac{\#\{ n \leq X : a(n) \equiv r \pmod{M} \}}{X}.
  \end{equation*}
An integral power series $A$ is called \textit{lacunary modulo $M$} if
\begin{equation*}
   \lim_{X \to \infty} \delta_0(A,M;X)=1,
\end{equation*}
which means that almost all the coefficients of $A$ are divisible by $M$. It turns out that such a result can also be proved for the $T_{\ell,k}(n)$ function. Specifically, we prove the following result in Section \ref{sec:lacunary}.

\begin{theorem}\label{thm1}
	Let $G(q)=\sum_{n \geq 0}T_{\ell,k}(n)q^n$, $k >0$ and $\ell >1$ be integers and $p^a$ be the largest power of a prime $p$ that divides $\ell$. If $p^{2a} \geq \ell$ and $k \leq p^{m+a}(1-p^{2s-2a})$, then for every positive integer $m$ and $s$ with $m > a > s \geq 0$, we have, 
 \begin{equation*}
     \lim_{X\to\infty} \delta_{0}(G,p^m;X)  = 1.
 \end{equation*}
\end{theorem}

\noindent The following is an easy corollary now.

\begin{corollary}
   Let $F(q)=\sum_{n \geq 0}T_{p,k}(n)q^n$. Then for every positive integer $m$ and for every prime $p$, we have
 \begin{equation*}
     \lim_{X\to\infty} \delta_{0}(F,p^m;X)  = 1.
 \end{equation*} 
\end{corollary}

Before proceeding to the proofs of our main theorems, we present a brief proof of the following result on the lacunarity of $T_2(9n+1)$.

\begin{theorem}\label{thm0.02365}
    The series $\sum_{n \geq 0} T_2(9n+1)q^n$ exhibits lacunarity modulo $6$, namely,
    \begin{equation*}
        \lim_{X \to \infty} \frac{\#\{ 0 \leq n \leq X : T_2(9n+1) \equiv 0 \pmod{6} \}}{X} = 1.
    \end{equation*}
\end{theorem}

To prove this, we recall the following classical result due to Landau \cite{lan}.

\begin{lemma}\label{l0.001}
    Let $r(n)$ and $s(n)$ be quadratic polynomials. Then
    \begin{equation*}
        \left( \sum_{n \in \mathbb{Z}} q^{r(n)} \right) \left( \sum_{n \in \mathbb{Z}} q^{s(n)} \right)
    \end{equation*}
    is lacunary modulo $2$.
\end{lemma}

\begin{proof}[Proof of Theorem \ref{thm0.02365}]
From \eqref{e50.1}, we have 
\begin{equation*}
    \sum_{n \geq 0} T_{2}(9n+1)q^n \equiv 3f_1f_2 \pmod{6}.
\end{equation*}
Using \eqref{e2.0.3.4} together with its magnified version (where $q \to q^2$) under modulo $2$ , we apply Lemma \ref{l0.001} to complete the proof of the theorem.
\end{proof}

The rest of the paper is organized as follows: In Section \ref{sec:prelim} we recall some results which we will use in our proofs, Sections \ref{sec:pf-1} -- \ref{sec:lacunary} contains the proofs of our results, we end the paper with some concluding remarks in Section \ref{sec:conc}.

\section{Preliminaries}\label{sec:prelim}

In this section we collect some results from elementary $q$-series analysis and the theory of modular forms which are useful in proving our results.

\subsection{Elementary Results}

We need Euler's pentagonal number theorem \cite[Eq. (1.3.18)]{Spirit}
\begin{equation}\label{e2.0.3.4}
    f_1=\sum_{n \in \mathbb{Z}}(-1)^nq^{n(3n-1)/2}.
\end{equation}
and Jacobi's Triple Product identity \cite[Eq. (1.3.24)]{Spirit}
\begin{equation}\label{e2.0.3.3}
    f_1^3=\sum_{n\geq 0}(-1)^n(2n+1)q^{n(n+1)/2}.
\end{equation}

Some known $3$-dissections are required (see for example \cite[Lemma 2.2]{DAS2025128913}).
\begin{lemma}
The following $3$-dissections holds
    \begin{align}
        \dfrac{f_1^2}{f_2} & = \dfrac{f_9^2}{f_{18}} -2q\dfrac{f_3f_{18}^2}{f_6f_9}, \label{e0.7}\\
        f_1^3 & = \dfrac{f_6f_9^6}{f_3f_{18^3}} -3q f_9^3 +4q^3 \dfrac{f_3^2f_{18}^6}{f_6^2f_9^3}, \label{e0.8}\\
        \dfrac{1}{f_1^3} & = a_3^2\dfrac{f_9^3}{f_{3}^{10}}+3qa_3\dfrac{f_9^6}{f_3^{11}}+9q^2\dfrac{f_9^9}{f_3^{12}},\label{e0.8.0}\\
     a_1& =a_3+6q\dfrac{f_9^3}{f_3}, \label{e0.7.0}
    \end{align}
where Borweins’ cubic theta function $a(q)$ is given by \[a_n = a(q^n) := \sum\limits_{j,k=-\infty}^{\infty}q^{n\cdot(j^2+jk+k^2)}.\]
\end{lemma}

We recall, partitions with even parts distinct wherein odd parts are unrestricted but even parts are distinct. Let $ped(n)$ denote the number of such partitions of $n$, then its generating function is given by
\begin{equation}\label{eq:gf-ped}
    \sum_{n\geq 0}ped(n)q^n=\frac{f_4}{f_1}.
\end{equation}
We need the following results related to this function.
\begin{lemma}\cite[Corollary 3.3]{andrews2010arithmetic}
    We have, for all $n\geq 0$
    \begin{align}
        ped(9n+7)& \equiv 0 \pmod{12}. \label{e2.6}
    \end{align}
\end{lemma}
\begin{lemma}\cite[Corollary 3.6]{andrews2010arithmetic}
    We have, for all $n\geq 0$
    \begin{align}
        ped\left( 3^{2\alpha +1 }n + \dfrac{17 \cdot 3^{2\alpha}-1}{8} \right) & \equiv 0 \pmod{6}, \label{e2.7}\\
        ped\left( 3^{2\alpha +2 }n + \dfrac{19 \cdot 3^{2\alpha+1}-1}{8} \right) & \equiv 0 \pmod{6}. \label{e2.8}
    \end{align}
\end{lemma}
\begin{lemma}\cite[Theorem 1.1]{Nath}
    We have, for all $n\geq 0$
    \begin{equation}\label{e3.2}
        ped(9n+7) \equiv ped\left( 9 \cdot 5^{2\alpha}n + \dfrac{57 \cdot 5^{2\alpha}-1}{8} \right) \pmod{24}.
    \end{equation}
\end{lemma}
\noindent We also need the following result from \cite[Lemma 7]{NadjiAhmia}.
\begin{lemma}
We have
\begin{align}
        \sum_{n \geq 0}T_{2}(3n+1)q^n & = 3\frac{f_2^4f_3^5}{f_1^8f_6}, \label{e0.2}\\
        \sum_{n \geq 0}T_{2}(3n+2)q^n & = 6\frac{f_2^3f_3^2f_6^2}{f_1^7}. \label{e0.2.1}
    \end{align}
\end{lemma}

The following result of Newman will play a crucial role in the proof of Theorem \ref{t2}, therefore we shall quote it as a lemma. Following the notations of Newman's paper, we shall let $p$, $q$ denote distinct primes, let $r, s \neq 0, r \not \equiv s \pmod{2}$. Set
\begin{equation}\label{e7}
    \phi(\tau) = \prod_{n=1}^{\infty}(1-x^n)^r(1-x^{nq})^s = \sum_{n=0}^{\infty}a(n)x^n,
\end{equation}
$\epsilon = \dfrac{1}{2}(r+s), t=(r+sq)/24, \Delta= t(p^2-1), \theta = (-1)^{\frac{1}{2}-\epsilon}2q^s$, then we have the following result.
\begin{lemma}\label{l2}\textup{\cite[Theorem 3]{16}}
    With the notations defined as above, the coefficients $c(n)$ of $\phi(\tau)$ satisfy
    \begin{equation*}
        a(np^2+\Delta)-\gamma(n)a(n) +p^{2\epsilon-2}a\left( \dfrac{n-\Delta}{p^2} \right) = 0,
    \end{equation*}
where 
\begin{equation*}
    \gamma(n) = p^{2\epsilon - 2}\alpha-\left( \dfrac{\theta}{p} \right)_L p^{\epsilon-3/2}\left( \dfrac{n-\Delta}{p} \right)_L,
\end{equation*}
and $\alpha$ is a constant.
\end{lemma}

\subsection{Preliminaries from Modular Forms}
We recall some basic facts and results from the theory of modular forms that will be useful in the proof of some of our results.

For a positive integer $N$, we will assume that:
\begin{align*}
\textup{SL}_2(\mathbb{Z}) & :=\left\{\begin{bmatrix}
a  &  b \\
c  &  d      
\end{bmatrix}: a, b, c, d \in \mathbb{Z}, ad-bc=1
\right\},\\
\Gamma_{0}(N) &:=\left\{
\begin{bmatrix}
a  &  b \\
c  &  d      
\end{bmatrix} \in \Gamma : c\equiv~0\pmod N \right\},\\
\Gamma_{1}(N) &:=\left\{
\begin{bmatrix}
a  &  b \\
c  &  d      
\end{bmatrix} \in \Gamma : a\equiv d \equiv 1\pmod N \right\},\\
\Gamma(N) &:=\left\{
\begin{bmatrix}
a  &  b \\
c  &  d      
\end{bmatrix} \in \textup{SL}_2(\mathbb{Z}) : a \equiv d \equiv 1 \pmod{N}, \text{and} \hspace{2mm} b\equiv c \equiv 0 \pmod N \right\}.
\end{align*}
A subgroup $\Gamma$ of $\textup{SL}_2(\mathbb{Z})$ is called a congruence subgroup if $\Gamma(N) \subseteq \Gamma$ for some $N$. The smallest $N$ such that $\Gamma(N) \subseteq \Gamma$ is called the level of $\Gamma$. For example, $\Gamma_{0}(N)$ and $\Gamma_{1}(N)$ are congruence subgroups of level $N$.\\

Let $\mathbb{H}:= \{ z \in \mathbb{C} : \Im(z) > 0 \}$ be the upper half of the complex plane. The group 

$$\textup{GL}_2^{+}(\mathbb{R}) = \left\{
\begin{bmatrix}
a  &  b \\
c  &  d      
\end{bmatrix} : a,b,c,d \in \mathbb{R} \hspace{2mm} \text{and} \hspace{2mm} ad-bc>0 \right\}$$ acts on $\mathbb{H}$ by $\begin{bmatrix}
a  &  b \\
c  &  d      
\end{bmatrix} z = \dfrac{az +b}{cz+d}$. We identify $\infty$ with $\dfrac{1}{0}$ and define $ \begin{bmatrix}
a  &  b \\
c  &  d      
\end{bmatrix} \dfrac{r}{s} = \dfrac{ar +bs}{cr+ds}$, where $\dfrac{r}{s} \in \mathbb{Q} \cup \{\infty\}$. This gives an action of $\textup{GL}_2^{+}(\mathbb{R})$ on the extended upper half-plane $\mathbb{H}^{\star} = \mathbb{H} \cup \mathbb{Q} \cup \{\infty\}$. Suppose that $\Gamma$ is a congruence subgroup of $\textup{SL}_2(\mathbb{Z})$. A cusp of $\Gamma$ is an equivalence class in $\mathbb{P}^{1}=\mathbb{Q} \cup \{ \infty\}$ under the action of $\Gamma$.\\

The group $\textup{GL}_2^{+}(\mathbb{R})$ also acts on the functions $f : \mathbb{H} \to \mathbb{C}$. In particular, suppose that $\gamma = \begin{bmatrix}
a  &  b \\
c  &  d      
\end{bmatrix} \in \textup{GL}_2^{+}(\mathbb{R}) $. If $f(z)$ is a meromorphic function on $\mathbb{H}$ and $\ell$ is an integer, then define the slash operator $\mid_{\ell}$ by 
$(f\mid_{\ell}\gamma )(z) := (det \gamma)^{\ell/2}(cz+d)^{-\ell}f(\gamma z)$.\\

\begin{definition}
Let $\gamma$ be a congruence subgroup of level N. A holomorphic function $f : \mathbb{H} \to \mathbb{C}$ is called a modular form with integer weight on $\Gamma$ if the following hold:
\begin{enumerate}
    \item We have
    \begin{align*}
        f\left(\frac{az+b}{cz+d}\right)=(cz+d)^{\ell}f(z)
    \end{align*}
    for all $z\in\mathbb{H}$ and all $\begin{bmatrix}
        a&b\\
        c&d
    \end{bmatrix}\in\Gamma$.
    \item If $\gamma\in\textup{SL}_2(\mathbb{Z})$, then $(f\mid_{\ell}\gamma )(z)$ has a Fourier expansion of the form
    \begin{align*}
        (f\mid_{\ell}\gamma )(z)=\sum_{n=0}^{\infty}a_{\gamma}(n)q^n_{N},
    \end{align*}
    where $q^n_{N}=e^{\frac{2\pi iz}{N}}$.
\end{enumerate}
\end{definition}

For a positive integer $\ell$, let $M_\ell(\Gamma_1(N)))$ denote the complex vector space of modular forms of weight $\ell$ with respect to $\Gamma_1(N)$. 
\begin{definition}\cite[Definition 1.15]{ono2004}
If $\chi$ is a Dirichlet character modulo $N$, then a modular form $f\in M_\ell(\Gamma_1(N))$ has Nebentypus character $\chi$ if 
$f\left( \frac{az+b}{cz+d}\right)=\chi(d)(cz+d)^{\ell}f(z)$ for all $z\in \mathbb{H}$ and all $\begin{bmatrix}
a  &  b \\
c  &  d      
\end{bmatrix}\in \Gamma_0(N)$. The space of such modular forms is denoted by $M_\ell(\Gamma_0(N),\chi)$.
\end{definition}

Finally, we require the following results to prove Theorem \ref{thm1}.
\begin{theorem}\cite[Theorem 1.64]{ono2004}\label{thm_ono1} Suppose that $f(z)=\displaystyle\prod_{\delta\mid N}\eta(\delta z)^{r_\delta}$ 
		is an eta-quotient such that $\ell=\displaystyle\dfrac{1}{2}\sum_{\delta\mid N}r_{\delta}\in \mathbb{Z}$, $\sum_{\delta\mid N} \delta r_{\delta}\equiv 0 \pmod{24}$, and  $\sum_{\delta\mid N} \dfrac{N}{\delta}r_{\delta}\equiv 0 \pmod{24}$.
		Then, 
		$
		f\left( \dfrac{az+b}{cz+d}\right)=\chi(d)(cz+d)^{\ell}f(z)
		$
		for every  $\begin{bmatrix}
			a  &  b \\
			c  &  d      
		\end{bmatrix} \in \Gamma_0(N)$. Here 
\(
		    \chi(d):=\left(\dfrac{(-1)^{\ell} \prod_{\delta\mid N}\delta^{r_{\delta}}}{d}\right)_L.\label{chi}\)
	\end{theorem}

\noindent If the eta-quotient $f(z)$ satisfies the conditions of Theorem \ref{thm_ono1} and  is holomorphic at all of the cusps of $\Gamma_0(N)$, then $f\in M_{\ell}(\Gamma_0(N), \chi)$. To determine  the orders of an eta-quotient at each cusp is the following.
	\begin{theorem}\cite[Theorem 1.65]{ono2004}\label{thm_ono1.1}
		Let $c, d,$ and $N$ be positive integers with $d\mid N$ and $\gcd(c, d)=1$. If $f(z)$ is an eta-quotient satisfying the conditions of Theorem~\ref{thm_ono1} for $N$, then the order of vanishing of $f(z)$ at the cusp $\dfrac{c}{d}$ 
		is $\dfrac{N}{24}\sum_{\delta\mid N}\dfrac{\gcd(d,\delta)^2r_{\delta}}{\gcd(d,\frac{N}{d})d\delta}.$
	\end{theorem}

We state the following deep result due to Serre, which is useful to prove  Theorem \ref{thm1}.
	\begin{theorem}\cite[Theorem~2.65]{ono2004}\label{serre}
		Let $k, m$  be positive integers. If  $f(z)\in M_{k}(\Gamma_0(N), \chi(\bullet))$ has the Fourier expansion $f(z)=\sum_{n \geq 0}c(n)q^n\in \mathbb{Z}[[q]],$
		then there is a constant $\alpha>0$  such that
		$$
		\# \left\{n\leq X: c(n)\not\equiv 0 \pmod{m} \right\}= \mathcal{O}\left(\dfrac{X}{\log^{\alpha}{}X}\right).
		$$
	\end{theorem}

As usual, we denote $M_k(SL_2(\mathbb{Z}))$ (resp. $S_k(SL_2(\mathbb{Z}))$) is the complex vector space of weight $k$ holomorphic modular (resp. cusp) forms with respect to $SL_2(\mathbb{Z})$.

\begin{definition}
Let $m$ be a positive integer and $f(z) = \sum_{n \geq 0} a(n)q^n \in M_{k}(\Gamma_0(N),\chi)$, where $\chi$ is a Dirichlet character modulo $N$. Then the action of Hecke operator $T_m$ on $f(z)$ is defined by 
\begin{align*}
f(z)|T_m := \sum_{n \geq 0} \left(\sum_{d\mid \gcd(n,m)}\chi(d)d^{k-1}a\left(\frac{nm}{d^2}\right)\right)q^n.
\end{align*}
In particular, if $m=p$ is prime, we have 
\begin{align}\label{hecke1}
f(z)|T_p := \sum_{n \geq 0} \left(a(pn)+\chi(p)p^{k-1}a\left(\frac{n}{p}\right)\right)q^n.
\end{align}
We note that $a(n)=0$ unless $n$ is a nonnegative integer.
\end{definition}
\begin{definition}\label{hecke2}
	A modular form $f(z)=\sum_{n \geq 0}a(n)q^n \in M_{k}(\Gamma_0(N),\chi)$ is called a Hecke eigenform if for every $m\geq2$ there exists a complex number $\lambda(m)$ for which 
	\begin{align}\label{hecke3}
	f(z)|T_m = \lambda(m)f(z).
	\end{align}
\end{definition}

In the special case when $\chi$ is the trivial character (i.e., $\chi(n) = 1$ for all $(n,N) = 1$), the Hecke operator simplifies to
\begin{align*}
f(z)|T_m := \sum_{n \geq 0} \left(\sum_{d\mid \gcd(n,m)}d^{k-1}a\left(\frac{nm}{d^2}\right)\right)q^n.
\end{align*}

The discriminant modular form $\Delta(z) = qf_1^{24} = \sum_{n \geq 1} \tau(n)q^n$ belongs to the space $S_{12}(SL_2(\mathbb{Z})) = S_{12}(\Gamma_0(1))$. Since the level is $N = 1$, the only Dirichlet character modulo $1$ is the trivial character. Therefore, the associated character $\chi$ is trivial. Here, the function $\tau(n)$ denotes the Ramanujan's tau function. It is known that $\Delta(z)$ is an eigenform for all Hecke operators. In other words,
\begin{equation*}
    \Delta(z) \mid T_{p} = \tau(p)\Delta(z)
\end{equation*}
for any prime $p$. The Ramanujan's tau function satisfies the following properties:
\begin{equation}\label{tau1}
    \tau(mn) = \tau(m)\tau(n), \quad \text{if} \quad \gcd(m,n)=1,
\end{equation}
\begin{equation}\label{tau2}
    \tau(p^{l}) = \tau(p)\tau(p^{l-1}) - p^{11}\tau(p^{l-2}),
\end{equation}
where $p$ is prime and $l \geq 2$.

\section{Proof of Theorem \ref{c1.4}}\label{sec:pf-1}
Considering \eqref{e0.1} and \eqref{e0.2} with $p=2$, $k=3$ and $l = 1$, we observe that
\begin{equation}\label{e0.3}
    \sum_{n\geq 0}T_{2}(3n+1)q^n \equiv 3 \dfrac{f_3^5}{f_6} \pmod{24}.
\end{equation}
Collecting the terms of the form $q^{3n+j}$ for $j=0,1,2$ from both sides of the equation \eqref{e0.3}, we get
\begin{align}
    \sum_{n\geq 0}T_{2}(9n+1)q^n & \equiv 3 \dfrac{f_1^5}{f_2} \pmod{24}, \label{e0.4}\\
    T_{2}(9n+4) &\equiv 0 \pmod{24}, \label{e0.5}\\
    T_{2}(9n+7) &\equiv 0 \pmod{24}. \label{e0.6}
\end{align}

Substituting \eqref{e0.7}, \eqref{e0.8} into \eqref{e0.4}, we obtain
\begin{equation*}
    \sum_{n\geq 0}T_{2}(9n+1)q^n \equiv 3\dfrac{f_6f_9^8}{f_3f_{18}^4}+9q \dfrac{f_9^5}{f_{18}} + 18q^2 \dfrac{f_3f_9^2f_{18}^2}{f_6} + 12q^3\dfrac{f_3^2f_{18}^5}{f_6^2f_9} \pmod{24}.
\end{equation*}
If we extract the terms involving $q^{3n+1}$ from both sides of the above equation, divide by $q$ and then replace $q^3$ by $q$, we get
\begin{equation}\label{e1.0}
    \sum_{n\geq 0}T_{2}(27n+10)q^n \equiv 9\dfrac{f_3^5}{f_{6}} \pmod{24}.
\end{equation}
Collecting the terms containing $q^{3n+j}$ for $j=0,1,2$ from both sides of \eqref{e1.0}, we get
\begin{align}
    \sum_{n\geq 0}T_{2}(81n+10)q^n & \equiv 9 \dfrac{f_1^5}{f_2} \pmod{24}, \label{e1.1}\\
    T_{2}(81n+37) &\equiv 0 \pmod{24}, \label{e1.2}\\
    T_{2}(81n+64) &\equiv 0 \pmod{24}. \label{e1.3}
\end{align}
Again by substituting \eqref{e0.7}, \eqref{e0.8} into \eqref{e1.1} and extracting the powers of the form $q^{3n+1}$ from both sides of the resulting equation, we find that
\begin{equation}\label{e1.4}
    \sum_{n\geq 0} T_{2}(243n+91)q^n \equiv 3 \dfrac{f_3^5}{f_6} \pmod{24}.
\end{equation}
From \eqref{e0.3}, \eqref{e1.0} and \eqref{e1.4}, we deduce that
\begin{align}
    3T_{2}(3n+1) &\equiv T_{2}(27n+10) \pmod{24} \label{e1.5},\\
    T_{2}(3n+1)  &\equiv T_{2}(243n+91) \pmod{24} \label{e1.6}.
\end{align}
Utilizing both \eqref{e1.5} and \eqref{e1.6} and by mathematical induction on $\alpha \geq 0$, we arrive at 
\begin{align}
    3T_{2}(3n+1) &\equiv T_{2}\left( 3^{4\alpha+3} + \sum_{i=0}^{2\alpha+1}3^{2i} \right) \pmod{24} \label{e1.7},\\
    T_{2}(3n+1)  &\equiv T_{2}\left( 3^{4\alpha+1} + \sum_{i=0}^{2\alpha}3^{2i} \right) \pmod{24} \label{e1.8}.
\end{align}
Using \eqref{e1.8}, \eqref{e0.5} and \eqref{e0.6}, we obtain \eqref{c0.1.4} and \eqref{c1.1.4}, respectively. Similarly, using \eqref{e1.7}, \eqref{e1.2} and \eqref{e1.3}, we obtain \eqref{c2.1.4} and \eqref{c3.1.4}, respectively.

\section{Proof of Theorem \ref{c1.4.1}}\label{sec:pf-2}
\begin{proof}[Proof of Theorem \ref{c1.4.1}]
We have
\begin{equation}\label{e1.9}
    \sum_{n\geq 0}T_{4}(n)q^n = \dfrac{f_4^3}{f_1^3}.
\end{equation}
Using \eqref{e0.1} in \eqref{e1.9} with $p=3$, $k=1$, we observe that
\begin{equation}\label{e2.0}
    \sum_{n\geq 0}T_{4}(n)q^n \equiv \dfrac{f_{12}}{f_3} \pmod{3}.
\end{equation}
Collecting the terms containing $q^{3n}$ from both sides of \eqref{e2.0}, we get
\begin{equation*}
    \sum_{n\geq 0}T_{4}(3n)q^n  \equiv \dfrac{f_4}{f_1} \pmod{3}, \label{e2.5}
\end{equation*}
Employing \eqref{eq:gf-ped}, \eqref{e2.6}, \eqref{e2.7}, \eqref{e2.8},  and \eqref{e3.2} into \eqref{e2.5}, we obtain \eqref{e3.0}, \eqref{e3.1} and \eqref{e2.9} respectively.
\end{proof}

As a consequence of Theorem \ref{c1.4.1}, we have the following corollary.
\begin{corollary}
    The function $T_{4}(n)$ is divisible by $3$ at least $13/648$ of the time.
\end{corollary}
\begin{proof}
The arithmetic sequences $3^{2\alpha +2 }n + \dfrac{17 \times 3^{2\alpha+1}-3}{8}$, $3^{2\alpha +3 }n + \dfrac{19 \times 3^{2\alpha+2}-3}{8}$ and $27 \cdot 5^{2\alpha}n + \dfrac{171 \cdot 5^{2\alpha}-3}{8}$ (for $\alpha \geq 1$), on which $T_{4}(\cdot)$ is $0$ modulo $3$, do not intersect. These sequences account for
\begin{equation*}
    \sum_{j=4}^{\infty}\dfrac{1}{3^j} + \dfrac{1}{27}\sum_{j=1}^{\infty}\dfrac{1}{5^{2j}} = \dfrac{13}{648}
\end{equation*}
of all positive integers.
\end{proof}

\section{Proof of Theorem \ref{t0.1.0.0}}\label{sec:thmn2}

From equations \eqref{e2.0.3.3} and \eqref{e0.2}, we have
    \begin{equation*}
        \sum_{n \geq 0} T_2(9n+1)q^n \equiv 3\sum_{j \geq 0}(-1)^n(2j+1)q^{j(j+1)/2} \pmod{6}.
    \end{equation*}
Therefore, if we wish to consider values of the form $T_2(9(pn+r)+1)$, then we need to know whether we can write $pn+r = j(j+1)/2$ for some nonnegative integer $j$. If we can show that no such representations exist, then the theorem is proved. Note that, if a representation of the form $pn+r = j(j+1)/2$ exists, then $r \equiv j(j+1)/2 \pmod{p}$. We complete the square to obtain $8r+1 \equiv (2j+1)^2 \pmod{p}$, which is not possible because we have assumed that $\left( \dfrac{8r+1}{p} \right)_L = -1$. Therefore, no such representation is possible, and this completes the proof.

\section{Proof of Theorem \ref{t0.0.1}}\label{s4}
We have the following, from \eqref{e0.2} and \eqref{e0.1}
\begin{equation}\label{e50.1}
    \sum_{n\geq0}T_{2}(9n+1)q^n \equiv 3f_1f_2 \pmod{6}.
\end{equation}
Combining \eqref{e2.0.3.4} and \eqref{e50.1}, we get
\begin{equation}\label{e50.3}
    \sum_{n\geq 0} T_{2}(9n+1)q^{24n+3} \equiv 3\sum_{m \in \mathbb{Z}}\sum_{k \in \mathbb{Z}} (-1)^{m+k}q^{{(6m-1)}^2+2(6k-1)^2} \pmod{6}.
\end{equation}
From \eqref{e50.3}, we find that if $24n+3$ is not of the form ${(6m-1)}^2+2(6k-1)^2$, then
\begin{equation*}
    T_2(9n+1) \equiv 0 \pmod{6}.
\end{equation*}
Let $\nu_p(N)$ denote the exponent of the highest power of $p$ dividing $N$. If $N$ is of the form $x^2+2y^2$, then due to the fact that $\left(\frac{-2}{p}\right)_L=-1$, it forces $p$ to divide both $x$ and $y$ and hence $\nu_p(N)$ must be even. Note that if $p \nmid n$, then $n$ can be written as $pk+j$ for some integer $k$ and $1\leq j\leq p-1$. Therefore, 
\begin{equation*}
    \nu_p\left(24\left( p^{2\alpha+1}(pk+j)+ \dfrac{p^{2\alpha+2}-1}{8} \right)+3 \right)=\nu_p\left(24\cdot p^{2\alpha+2}k+24\cdot j\cdot p^{2\alpha+1}+3\cdot p^{2\alpha+2} \right)=2\alpha+1.
\end{equation*}
So, $ 24 \left( p^{2\alpha+1}n+ \dfrac{p^{2\alpha+2}-1}{8} \right)+3$ is not of the form $x^2+2y^2$ when $p\nmid n$ and hence \eqref{e50.0} holds.

\section{Proof of theorem \ref{t2}}\label{s5}
From \eqref{e0.2.1} and using \eqref{e0.1}, we have that
    \begin{equation}\label{e10}
        \sum_{n \geq 0}T_{2}(3n+2)q^n = 6\dfrac{f_2^3f_3^2f_6^2}{f_1^7} \equiv 6 \dfrac{f_3^6}{f_1} \equiv 6  \sum_{n \geq 0}a(n)q^n \pmod{12},
    \end{equation}
    where $\dfrac{f_3^6}{f_1} = \sum_{n \geq 0}a(n)q^n$. Since there is a factor of 6 on the right-hand side of \eqref{e10}, it is sufficient to prove the corresponding congruences for $a(n)$ modulo 2.
    
    Putting $r=-1, q=3$ and $s=6$ in $\eqref{e7}$, we get $\epsilon = \frac{5}{2}, \Delta=\frac{17}{24}(p^2-1), \theta=2\cdot3^6$. Therefore, by Lemma $\ref{l2}$, for any $n\geq 0$, we have
    \begin{equation}\label{e11}
        a\left( p^2n +\dfrac{17}{24}(p^2-1) \right) - \gamma(n)a(n) + p^3a\left( \dfrac{1}{p^2}\left( n-\dfrac{17}{24}(p^2-1) \right) \right)=0,
    \end{equation}
    where
    \begin{equation}\label{e12}
        \gamma(n)=p^3\alpha-p\left(\dfrac{2}{p}\right)_L\left(\dfrac{n-\frac{17}{24}(p^2-1)}{p}\right)_L,
    \end{equation}
    and $\alpha$ is a constant integer. Setting $n=0$ in $\eqref{e11}$ and using the fact that $a(0) = 1$ and \[a\left(\dfrac{\frac{-17}{24}(p^2-1)}{p^2}\right)=0,\] we obtain
    \begin{equation}\label{e13}
        a\left( \dfrac{17}{24}(p^2-1) \right) = \gamma(0).
    \end{equation}
Setting $n=0$ in $\eqref{e12}$ and using $\eqref{e13}$, we obtain
\begin{equation}\label{e14}
    p^3\alpha = a\left( \frac{17}{24}(p^2-1) \right)+p\left(
    \dfrac{2}{p}\right)_L \left(\dfrac{\frac{-17}{24}(p^2-1)}{p}\right)_L := \omega(p).
    \end{equation}
Now rewriting $\eqref{e11}$, by referring to $\eqref{e12}$ and $\eqref{e14}$, we obtain
\begin{align}
        a\left( p^2n +\dfrac{17}{24}(p^2-1) \right) &= \left( \omega(p)-p\left(\dfrac{2}{p}\right)_L\left( \dfrac{n-\frac{17}{24}(p^2-1)}{p} \right)_L \right)a(n) \nonumber\\
        &\quad- p^3a\left( \dfrac{1}{p^2}\left( n-\dfrac{17}{24}(p^2-1) \right) \right).\label{e15}
    \end{align}
Now, replacing $n$ by $pn + \dfrac{17}{24}(p^2-1)$ in $\eqref{e15}$, we obtain
\begin{equation}\label{e16}
    a\left( p^3n+\dfrac{17}{24}(p^4-1) \right) = \omega(p)a\left( pn + \dfrac{17}{24}(p^2-1) \right) - p^3 a(n/p).
\end{equation}

From equation $\eqref{e11}$, we can see that
\begin{equation}\label{e28}
        a\left( p^2n +\dfrac{17}{24}(p^2-1) \right) - \gamma(n)a(n) + a\left( \dfrac{1}{p^2}\left( n-\frac{17}{24}(p^2-1) \right) \right) \equiv 0 \pmod{2},
    \end{equation}
    where
    \begin{equation}\label{e29}
        \gamma(n) \equiv \omega(p) + \left(\dfrac{n-\frac{17}{24}(p^2-1)}{p}\right)_L  \pmod2.
    \end{equation}
Setting $n=0$ in $\eqref{e28}$ and using the fact that $a(0)=1$ and $a\left( \dfrac{\frac{-17}{24}(p^2-1)}{p^2} \right) = 0$, we arrive at
\begin{equation}\label{e30}
    a\left( \dfrac{17}{24}(p^2-1) \right) \equiv \gamma(0) \pmod{2}.
\end{equation}
Setting $n=0$ in $\eqref{e29}$ yields
\begin{equation}\label{e31}
    \gamma(0) \equiv \omega(p) + 1 \pmod2, \qquad \text{for} \quad p \neq 17,
\end{equation}
and
\begin{equation}\label{e31.1}
    \gamma(0) \equiv \omega(p) \pmod2, \qquad \text{for} \quad p = 17.
\end{equation}
Combining $\eqref{e30}$ and $\eqref{e31}$ yields
\begin{equation}\label{e32}
    a\left( \dfrac{17}{24}(p^2-1) \right) + 1 \equiv \omega(p) \pmod{2}, \quad p\neq 17,
\end{equation}
and combining $\eqref{e30}$ and $\eqref{e31.1}$ yields
\begin{equation}\label{e32:2}
    a\left( \dfrac{17}{24}(p^2-1) \right) \equiv \omega(p) \pmod{2}, \quad p=17.
\end{equation}
\underline{\textbf{Case - 1 :} $\omega(p) \equiv 0 \pmod{2}$}

From $\eqref{e16}$ we obtain
\begin{equation}\label{e17}
    a\left( p^3n+\dfrac{17}{24}(p^4-1) \right) \equiv p^3 a(n/p) \pmod{2}.
\end{equation}
Now, replacing $n$ by $pn$ in $\eqref{e17}$, we obtain
\begin{equation*}
    a\left(p^4n+\dfrac{17}{24}(p^4-1)\right) \equiv p^3 a(n) \equiv a(n) \pmod{2}.
\end{equation*}
Since $p^{4k}n+\dfrac{17}{24}(p^{4k}-1)=p^4\left(p^{4k-4}n+\dfrac{17}{24}(p^{{4k-4}}-1)\right)+\dfrac{17}{24}(p^4-1)$, using equation $\eqref{e17}$, we obtain that for every integer $k\geq 1$,
\begin{equation}\label{e19}
    a\left(p^{4k}n+\dfrac{17}{24}(p^{4k}-1)\right) \equiv a\left(p^{4k-4}n+\dfrac{17}{24}(p^{4k-4}-1)\right) \equiv a(n) \pmod{2}.
\end{equation}
Now if $p \nmid n$, then $\eqref{e17}$ yields
\begin{equation}\label{e20}
    a\left(p^3n +\dfrac{17}{24}(p^4-1)\right) \equiv 0 \pmod{2}.
\end{equation}
Replacing $n$ by $p^3+\dfrac{17}{24}(p^4-1)$ in $\eqref{e19}$ and using $\eqref{e20}$, we obtain, for $p \nmid n$, that
\begin{equation}\label{e21}
a\left( p^{4k+3}n + \dfrac{17}{24}(p^{4k+4}-1) \right) \equiv 0 \pmod{2}.    
\end{equation}
In particular, for $1 \leq j \leq p-1$, replacing $n$ by $pn+j$ we have from $\eqref{e21}$
\begin{equation}\label{e40}
    a\left( p^{4k+4}n + p^{4k+3}j + \dfrac{17}{24}\left( p^{4k+4}-1 \right)\right) \equiv 0 \pmod{2}.
\end{equation}
Congruence $\eqref{t3.1}$ follows from \eqref{e10} and $\eqref{e40}$. \\

Now, assume $p=17$, we have from $\eqref{e29}$ and $\eqref{e32:2}$, that
\begin{equation}\label{e51}
    \gamma(n) \equiv a\left( \dfrac{17}{24}(p^2-1) \right) + \left( \dfrac{n-\frac{17}{24}(p^2-1)}{p} \right)_L \pmod{2}.
\end{equation}
By $\eqref{e28}$ and $\eqref{e51}$
\begin{align}
    a\left( p^2n + \dfrac{17}{24}(p^2-1) \right) &\equiv \left( a\left( \dfrac{17}{24}(p^2-1) \right) + \left( \dfrac{n - \frac{17}{24}(p^2-1)}{p} \right)_L \right) a(n) \nonumber \\
    &\quad + a\left( \dfrac{1}{p^2} \left( n - \dfrac{17}{24}(p^2-1) \right) \right) \pmod{2}. \label{e52}
\end{align}
Replacing $n$ by $pn + \dfrac{17}{24}(p^2-1)$ in $\eqref{e52}$ yields
\begin{equation}\label{e53}
    a\left( p^3n + \dfrac{17}{24}(p^4-1) \right) \equiv a\left( \dfrac{17}{24}(p^2-1) \right)a\left( pn + \dfrac{17}{24}(p^2-1) \right) + a(n/p) \pmod{2}.
\end{equation}

For $p=17$, we have $\omega(p) \equiv 0 \pmod{2}$ so from $\eqref{e32:2}$, we have $a\left( \dfrac{17}{24}(p^2-1) \right) \equiv 0 \pmod{2}$. Therefore, replacing $n$ by $pn$ in $\eqref{e53}$ we find that
\begin{equation}\label{e54}
    a\left( p^4n + \dfrac{17}{24}(p^4-1) \right) \equiv a(n) \pmod{2}.
\end{equation}
By $\eqref{e54}$ and iteration, we deduce that for $n,k \geq 0$,
\begin{equation}\label{e55}
    a\left( p^{4k}n + \dfrac{17}{24}(p^{4k}-1) \right) \equiv a(n) \pmod{2}.
\end{equation}
Moreover, we can rewrite $\eqref{e52}$ as 
\begin{equation}\label{e56}
    a\left( p^2n + \dfrac{17}{24}(p^2-1) \right) \equiv a(n) \left( \dfrac{n-\frac{17}{24}(p^2-1)}{p} \right)_L + a\left( \dfrac{1}{p^2}\left( n-\dfrac{17}{24}(p^2-1) \right) \right) \pmod{2}.
\end{equation}
If $p\nmid (24n+17)$ and $n \not\equiv 0 \pmod{p}$, then $p \nmid \left( n -\dfrac{17}{24}(p^2-1) \right)$ and $\dfrac{n-\frac{17}{24}(p^2-1)}{p^2}$ is not an integer. Note that we require $n \not\equiv 0 \pmod{p}$, since otherwise $p \mid \left(n - \frac{17}{24}(p^2 - 1)\right)$, which is not desired. Therefore
\begin{equation}\label{e57}
    \left( \dfrac{n-\frac{17}{24}(p^2-1)}{p} \right)_L \equiv 1 \pmod{2},
\end{equation}
and
\begin{equation}\label{e58}
    a\left( \dfrac{n-\frac{17}{24}(p^2-1)}{p^2} \right) = 0.
\end{equation}
On account of \eqref{e56}, \eqref{e57}, and \eqref{e58}, we deduce that
\begin{equation}\label{e59}
    a\left( p^2n + \dfrac{17}{24}(p^2-1) \right) \equiv a(n) \pmod{2}.
\end{equation}
Replacing $n$ by $p^2n + \dfrac{17}{24}(p^2-1)$ in $\eqref{e55}$ and employing $\eqref{e59}$, we see that for $k\geq 0$,
\begin{equation}\label{e60}
    a\left( p^{4k+2}n + \dfrac{17}{24}(p^{4k+2}-1) \right) \equiv a(n) \pmod{2}.
\end{equation}
Equations $\eqref{e10}$ and $\eqref{e60}$ readily yield $\eqref{t3.3.1}$.\\

\underline{\textbf{Case - 2 :} $\omega(p) \equiv 1 \pmod{2}$}
In order to prove $(ii)$, we replace $n$ by $p^2n+\dfrac{17}{24}p(p^2-1)$ in $\eqref{e16}$.
\begin{align}
    a\left( p^5n + \dfrac{17}{24}(p^6-1) \right) & = a\left( p^3\left( p^2n+\dfrac{17}{24}p(p^2-1) \right) + \dfrac{17}{24}(p^4-1) \right) \nonumber\\
    &= \omega(p)a\left(p^3n+\dfrac{17}{24}(p^4-1)\right)-p^3a\left( pn + \dfrac{17}{24}(p^2-1) \right) \nonumber \\
    & = \left[ \omega^2(p) - p^3\right]a\left( pn+\dfrac{17}{24}(p^2-1) \right) -p^3\omega(p)a(n/p), \label{e22}
\end{align}
where the last equality follows from \eqref{e16}.

Now, as $\omega(p) \equiv 1 \pmod{2}$ and $p\geq5$ is an odd prime, we have $\omega^2(p)-p^3\equiv0 \pmod{2}$, and therefore $\eqref{e22}$ becomes
\begin{equation}\label{e23}
    a\left( p^5n + \dfrac{17}{24}(p^6-1) \right) \equiv a(n/p) \pmod{2}.
\end{equation}
Replacing $n$ by $pn$ in $\eqref{e23}$, we obtain
\begin{equation}\label{e24}
    a\left( p^6n +\dfrac{17}{24}(p^6-1) \right) \equiv a(n) \pmod{2}.
\end{equation}
Using equation $\eqref{e24}$ repeatedly, we see that for every integers $k\geq 1$,
\begin{equation}\label{e25}
    a\left( p^{6k}n + \dfrac{17}{24}(p^{6k}-1) \right) \equiv a(n) \pmod{2}.
\end{equation}
Observe that if $p\nmid n$, then $a(n/p)=0$. Thus $\eqref{e23}$ yields
\begin{equation}\label{e26}
    a\left( p^5n + \dfrac{17}{24}(p^6-1) \right) \equiv 0 \pmod{2}.
\end{equation}
Replacing $n$ by $p^5n+\dfrac{17}{24}(p^6-1)$ in $\eqref{e25}$ and using $\eqref{e26}$, we obtain, for $p \nmid n$, that
\begin{equation}\label{e27}
    a\left( p^{6k+5}n + \dfrac{17}{24}(p^{6k+6}-1) \right) \equiv 0 \pmod{2}.
\end{equation}
In particular, for $1 \leq j \leq p-1$, replacing $n$ by $pn+j$ we have from $\eqref{e27}$
\begin{equation}\label{e27.1}
    a\left( p^{6k+6}n + p^{6k+5}j + \dfrac{17}{24}\left( p^{6k+6}-1 \right)\right) \equiv 0 \pmod{2}.
\end{equation}
Congruence $\eqref{t2.2}$ readily follows from $\eqref{e10}$ and $\eqref{e27.1}$. \\

Again, from $\eqref{e29}$  and $\eqref{e32}$, we have
\begin{equation}\label{e33}
    \gamma(n) \equiv a\left( \dfrac{17}{24}(p^2-1)\right) + 1 + \left( \dfrac{n-\frac{17}{24}(p^2-1)}{p} \right)_L \pmod{2}.
\end{equation}
By $\eqref{e28}$ and $\eqref{e33}$,
\begin{align}
    a\left( p^2n + \dfrac{17}{24}(p^2-1) \right) &\equiv \left( a\left( \dfrac{17}{24}(p^2-1) \right) + 1 + \left( \dfrac{n - \frac{17}{24}(p^2-1)}{p} \right)_L \right) a(n) \nonumber \\
    &\quad + a\left( \dfrac{1}{p^2} \left( n - \dfrac{17}{24}(p^2-1) \right) \right) \pmod{2} \label{e34}.
\end{align}

Replacing $n$ by $pn + \dfrac{17}{24}(p^2-1)$ in $\eqref{e34}$ yields 
\begin{equation}\label{e35}
    a\left( p^3n + \dfrac{17}{24}(p^4-1) \right) \equiv \left( a\left( \dfrac{17}{24}(p^2-1) \right) + 1 \right)a\left( pn + \dfrac{17}{24}(p^2-1) \right) + a(n/p).
\end{equation}
Since $\omega(p) \equiv 1 \pmod{2}$ so from $\eqref{e32}$, we have $a\left( \dfrac{17}{24}(p^2-1) \right) \equiv 0 \pmod{2}$. Therefore, replacing $n$ by $pn$ in $\eqref{e35}$, we find that 
\begin{equation}\label{e42}
    a\left( p^{4}n + \dfrac{17}{24}(p^4-1) \right) \equiv a\left( p^2n+\dfrac{17}{24}(p^2-1) \right) + a(n) \pmod{2}.
\end{equation}
Replacing $n$ by $p^2n + \dfrac{17}{24}(p^2-1)$ in $\eqref{e42}$ yields
\begin{equation}\label{e43}
    a\left( p^6n + \dfrac{17}{24}(p^6-1) \right) \equiv a\left( p^4n + \dfrac{17}{24}(p^4-1) \right) + a\left( p^2n+ \dfrac{17}{24}(p^2-1) \right) \pmod{2}.
\end{equation}
Combining $\eqref{e42}$ and $\eqref{e43}$ yields
\begin{equation}\label{e44}
    a\left( p^6n + \dfrac{17}{24}(p^6-1) \right) \equiv a(n) \pmod{2}.
\end{equation}
By $\eqref{e44}$ and iteration, we deduce that for $n,k \geq 0$,
\begin{equation}\label{e45}
    a\left( p^{6k}n + \dfrac{17}{24}(p^{6k}-1) \right) \equiv a(n) \pmod{2}.\\
\end{equation}

Moreover, we can rewrite $\eqref{e34}$ as 
\begin{multline}
    \label{e46}
    a\left( p^2n + \dfrac{17}{24}(p^2-1) \right) \equiv \left(1 + \left( \dfrac{n-\frac{17}{24}(p^2-1)}{p} \right)_L \right) a(n)\\ + a\left( \dfrac{1}{p^2}\left( n-\dfrac{17}{24}(p^2-1) \right) \right) \pmod{2}.
\end{multline}

If $p\nmid (24n+17)$ and $n \not\equiv -17 \cdot 24^{-1} \pmod{p}$, then $p \nmid \left( n -\dfrac{17}{24}(p^2-1) \right)$ and $\dfrac{n-\frac{17}{24}(p^2-1)}{p^2}$ is not an integer. Note that we require $n \not\equiv -17 \cdot 24^{-1} \pmod{p}$, since otherwise $p \mid \left(n - \frac{17}{24}(p^2 - 1)\right)$, which is not desired. We note that equations \eqref{e57} and \eqref{e58} now follow. On account of \eqref{e46}, \eqref{e57}, and \eqref{e58}, we deduce that
\begin{equation}\label{e49}
    a\left( p^2n + \dfrac{17}{24}(p^2-1) \right) \equiv 0 \pmod{2}.
\end{equation}
Replacing $n$ by $p^2n + \dfrac{17}{24}(p^2-1)$ in $\eqref{e45}$ and employing $\eqref{e49}$, we see that for $k\geq 0$,
\begin{equation}\label{e50}
    a\left( p^{6k+2}n + \dfrac{17}{24}(p^{6k+2}-1) \right) \equiv 0 \pmod{2}.
\end{equation}
Congruence $\eqref{t3.3}$ follows from $\eqref{e10}$ and $\eqref{e50}$.\\

\section{Proof of Theorem \ref{t4}}\label{s3}
In this section, we prove Theorem \ref{t4} using the theory of Hecke eigenforms.

\noindent Magnifying $q \to q^8$ and multiplying $q$ on both sides of $\eqref{e1.0.0.0}$, we have
\begin{equation}\label{e4}
    \sum_{n \geq 0}T_{2}(n)q^{8n+1}=q\dfrac{f_{16}^3}{f_8^3} \equiv \Delta(z) = \sum_{n \geq 1}\tau(n)q^n \pmod{2}.
\end{equation}
Hence, from $\eqref{e2.0.3.3}$, we have
\begin{equation*}
    \Delta(z)=qf_1^{24} \equiv qf_8^3 \equiv \sum_{n \geq 0}q^{(2n+1)^2} \pmod{2}.
\end{equation*}
Therefore, we have
\begin{equation}\label{e5}
    \sum_{n \geq 0}T_{2}(n)q^{8n+1}\equiv \sum_{n \geq 0}q^{(2n+1)^2} \pmod{2}.
\end{equation}
If $s\equiv1 \pmod{8}$ and $\left(\dfrac{s}{p}\right)_L = -1$, then, for any $n\geq 0$, $8np+s$ cannot be an odd square. This implies that the coefficients of $q^{8np+s}$ in the left-hand side of $\eqref{e5}$ must be even. It follows that
\begin{equation}\label{e6}
    T_{2}\left( pn + \dfrac{s-1}{8} \right) \equiv 0 \pmod{2}.
\end{equation}
This proves $\eqref{e12.0.1}$.\\

Since $\tau(p) \equiv 0 \pmod{2}$ and $\Delta(n)$ is a Hecke eigenform, we have
\begin{equation*}
    \Delta(z)\mid T_{p} = \tau(p)\Delta(z) \equiv 0 \pmod{2}.
\end{equation*}
By $\eqref{hecke1}$ and $\eqref{e4}$, we get
\begin{equation*}
    \sum_{n \geq 0}T_{2}(n) q^{8n+1} \mid T_{p} \equiv \sum_{n \geq 0}\left( T_{2}\left( \dfrac{pn-1}{8} \right) + T_{2}\left( \dfrac{n/p-1}{8}\right) \right)q^n \equiv 0 \pmod{2}.
\end{equation*}
If we write $m=\dfrac{n/p-1}{8}\in N $, then $\dfrac{pn-1}{8} = p^2m + \dfrac{p^2-1}{8}$. So, we have
\begin{equation*}
    T_{2}\left( p^2m + \dfrac{p^2-1}{8} \right) + T_{2}(m) \equiv 0 \pmod{2}.
\end{equation*}
That is, 
\begin{equation*}
    T_{2}\left( p^2m + \dfrac{p^2-1}{8} \right) \equiv T_{2}(m) \pmod{2}.
\end{equation*}
By induction, for $k\geq 1$ we find that
\begin{equation}\label{ef1}
    T_{2}\left( p^{2k}m + \dfrac{p^{2k}-1}{8} \right) \equiv T_{2}(m) \pmod{2}.
\end{equation}
Replacing $m$ by $pn + \dfrac{s-1}{8}$ in \eqref{ef1}, and considering $\eqref{e6}$, we get
\begin{equation*}
    T_{2}\left( p^{2k+1}n+\dfrac{sp^{2k}-1}{8} \right) \equiv T_{2}\left( pn + \dfrac{s-1}{8} \right) \equiv 0 \pmod{2}.
\end{equation*}
This proves $\eqref{e12.0.2}$.\\

Since $\tau(p) \equiv 0 \pmod{2}$, $\eqref{tau2}$ gives
\begin{equation*}
    \tau(p^{2k+1}) \equiv 0 \pmod{2}.
\end{equation*}
Applying $\eqref{tau1}$, we have
\begin{equation*}
    \tau\left( p^{2k+1} \left( 8pn + r \right) \right) = \tau\left( p^{2k+1} \right)\tau\left( 8pn+r \right) \equiv 0 \pmod{2}.
\end{equation*}
It follows from the above relation and $\eqref{e4}$ that
\begin{equation*}
    T_{2}\left( p^{2k+2}n + \dfrac{rp^{2k+1}-1}{8} \right) \equiv 0 \pmod{2}.
\end{equation*}
This proves $\eqref{e12.0.3}$.

\section{Proof of Theorem \ref{thm1.00}}\label{sec:thmn1}

We use the theory of Hecke eigenforms to prove Theorems \ref{thm1.00}.\\

We have from \eqref{e50.1} that
		\begin{align*}
		\sum_{n \geq 0}T_2(9n+1)q^n \equiv 3f_1^3 \pmod{6}.
		\end{align*}
This gives
\begin{align*}
\sum_{n \geq 0}T_2(9n+1)q^{8n+1} \equiv 3\eta(4z)^6 \pmod{6}.
\end{align*}	
Let $\eta(4z)^6 = \sum_{n \geq 0} a(n)q^n.$ Then $a(n) = 0$ if $n\not\equiv 1\pmod{8}$ and for all $n\geq 0$, 
\begin{align}\label{1_1}
T_2(9n+1) \equiv a(8n+1) \pmod 6.
\end{align}	
By Theorem \ref{thm_ono1}, we have $\eta(4z)^6 \in S_3(\Gamma_0(16))$. Since $\eta(4z)^6$ is a Hecke eigenform (see, for example \cite{martin}), \eqref{hecke1} and \eqref{hecke3} yield
\begin{align*}
	\eta(4z)^6|T_p = \sum_{n \geq 1} \left(a(pn) + \chi(p)p^2 a\left(\frac{n}{p}\right) \right)q^n = \lambda(p) \sum_{n \geq 1} a(n)q^n,
\end{align*}
which implies 
\begin{align}\label{1.1}
	a(pn) + \chi(p)p^2 a\left(\frac{n}{p}\right) = \lambda(p)a(n).
\end{align}
Putting $n=1$ and noting that $a(1)=1$, we readily obtain $a(p) = \lambda(p)$.
Since $a(p)=0$ for all $p \not\equiv 1 \pmod{8}$, we have $\lambda(p) = 0$.
From \eqref{1.1}, we obtain 
\begin{align}\label{new-1.1}
a(pn) + \chi(p)p^2 a\left(\frac{n}{p}\right) = 0.
\end{align}
From \eqref{new-1.1}, we derive that for all $n \geq 0$ and $p\nmid r$, 
\begin{align}\label{1.2}
a(p^2n + pr) = 0
\end{align}
and  
\begin{align}\label{1.3}
a(p^2n) = - \chi(p)p^2 a(n)\equiv a(n) \pmod{2}.
\end{align}
Substituting $n$ by $8n-pr+1$ in \eqref{1.2} and together with \eqref{1_1}, we find that
\begin{align}\label{1.4}
T_2\left(9p^2n + \frac{9p^2-1}{8}+ 9pr\frac{1-p^2}{8}\right) \equiv 0 \pmod 6.
\end{align}	
Substituting $n$ by $8n+1$ in \eqref{1.3} and using \eqref{1_1}, we obtain
\begin{align}\label{1.5}
T_2\left(9p^2n + \frac{9p^2-1}{8}\right) \equiv T_2(9n+1) \pmod 6.
\end{align}
Since $p \geq 5$ is prime, so $8\mid (1-p^2)$ and $\gcd \left(\frac{1-p^2}{8} , p\right) = 1$.  
Hence when $r$ runs over a residue system excluding the multiple of $p$, so does $\frac{1-p^2}{8}r$. Thus \eqref{1.4} can be rewritten as
\begin{align}\label{1.6}
T_2\left(9p^2n + \frac{9p^2-1}{8}+ 9pj\right) \equiv 0 \pmod 6,
\end{align}
where $p \nmid j$. 
\par Now, $p_i \geq 5$ are primes such that $p_i \not\equiv 1 \pmod 8$. Since 
\begin{align*}
9p_1^2\dots p_{k}^2n + \frac{9p_1^2\dots p_{k}^2-1}{8}=9p_1^2\left(p_2^2\dots p_{k}^2n + \frac{p_2^2\dots p_{k}^2-1}{8}\right)+\frac{9p_1^2-1}{8},
\end{align*}
using \eqref{1.5} repeatedly we obtain that
\begin{align*}
T_2\left(9p_1^2\dots p_{k}^2n + \frac{9p_1^2\dots p_{k}^2-1}{8}\right) \equiv T_2(9n+1) \pmod{6},
\end{align*}
which implies
\begin{equation}\label{1.7}
    T_2\left( p_1^2 \cdots p_k^2 n + \frac{p_1^2\cdots p_k^2 -1}{8} \right) \equiv T_2(n) \pmod{6}.
\end{equation}
Let $j\not\equiv 0\pmod{p_{k+1}}$. Then \eqref{1.6} and \eqref{1.7} yield
\begin{align*}
T_2\left(9p_1^2\dots p_{k+1}^2n + \frac{9p_1^2\dots p_{k}^2p_{k+1}(8j+p_{k+1})-1}{8}\right) \equiv 0 \pmod{6}.
\end{align*}
This completes the proof of the theorem.

\section{Proof of Theorem \ref{thm1}}\label{sec:lacunary}
 To prove Theorem \ref{thm1}, we need the following lemmas.
	\begin{lemma}\label{lem2}For any prime $p$, positive integers $ \ell \geq 2 $, $ k \geq 1 $, and $ a \geq 1 $, and for any positive integer $ m > 2 $, we have
		\begin{align*}
			 \dfrac{\eta(24\ell z)^k \eta(24z)^{p^{a+m}-k}}{\eta(24p^az)p^m} \equiv \sum_{n \geq 0}T_{\ell,k}(n)q^{24n+k(\ell-1)} \pmod {p^{m}}.
		\end{align*}
	\end{lemma}
	\begin{proof} 	
		Consider 
		\begin{align*}
			\mathcal{A}(z) 
			=\dfrac{\eta(24z)^{p^a}}{\eta(24p^a z)}.
		\end{align*}
		We have
		\begin{align*}
			\mathcal{A}^{p^{m}}(z) = \dfrac{\eta(24z)^{p^{m+a}}}{\eta(24p^a z)^{p^m}} \equiv 1 \pmod {p^{m+1}} \quad \text{(follows from \eqref{e0.1} and \eqref{eta 1}).}
		\end{align*}
		Define $\mathcal{B}_{\ell, p, m}(z)$ by
		$$\mathcal{B}_{\ell, p, m}(z)=	\dfrac{\eta(24\ell z)^k}{\eta(24z)^k}\mathcal{A}^{p^{m}}(z).$$
		\noindent Now, working modulo $p^{m+1}$, we have 
		\begin{align}\label{new-110}
			\mathcal{B}_{\ell, p, m}(z)	&=\dfrac{\eta(24\ell z)^k}{\eta(24z)^k}\dfrac{\eta(24z)^{p^{m+a}}}{\eta(24p^az)^{p^m}}\equiv \dfrac{\eta(24\ell z)^k}{\eta(24z)^k}= q^{k(\ell-1)}\dfrac{f_{24\ell}^k}{f_{24}^k} \pmod{p^{m+1}}.
		\end{align}
Combining \eqref{eq:gf-lk} and \eqref{new-110}, we obtain the required result.
	\end{proof}
	\begin{lemma}\label{lem1}
		Let $\ell$ be an integer greater than one, $p$ be a prime divisor of $\ell$, $a$ be the largest positive integer such that $p^a$ divides $\ell$, and $t$ be a positive integer such that $t$ divides $\ell$. If $p^{2a} \geq \ell$ and $k \leq p^{m+a}(1-p^{2s-2a})$, then for every positive integer $m$, we have
		\begin{align*}
			\mathcal{B}_{\ell, p, m}(z) \in M_{\frac{p^m(p^a-1)}{2}}\left(\Gamma_0(576 \cdot \ell), \chi_1(\bullet)\right),
		\end{align*}
	where the Nebentypus character $$\chi_1(\bullet)=\left(\dfrac{(-1)^{\frac{p^m(p^a-1)}{2}} (24\ell)^{k}\cdot (24)^{p^{m+a}-k} \cdot (24p^a)^{-p^m}}{\bullet}\right)_L.$$	
	\end{lemma}
\begin{proof}
    First we verify the first, second and third hypotheses of Theorem \ref{thm_ono1}. The weight of the eta-quotient is $\dfrac{p^m(p^a-1)}{2}$. Suppose the level of the eta-quotient $\mathcal{B}_{\ell, p, m}(z)$ is $24\ell u$, where $u$ is the smallest positive integer satisfying the following identity.
    \begin{equation*}
        24 \ell u \left( \dfrac{k}{24\ell}+ \dfrac{p^{m+a}-k}{24} - \dfrac{p^m}{24p^a} \right) \equiv 0 \pmod{24}.
    \end{equation*}
Equivalently, we have
\begin{equation}\label{e100}
    u\left( k(1-\ell)+\ell p^{m-a}(p^{2a}-1) \right) \equiv0 \pmod{24}.
\end{equation}
Hence, from \eqref{e100}, we conclude that the level of the eta-quotient $\mathcal{B}_{\ell, p, m}$ is $576 \ell$ for $m>a$. By Theorem \ref{thm_ono1.1}, the cusps of $\Gamma_0(576\ell)$ are given by $\dfrac{c}{d}$ where $d \mid 576\ell$ and $\text{gcd}(c,d)=1$. Now $\mathcal{B}_{\ell, p, m}$ is holomorphic at a cusp $\dfrac{c}{d}$ if and only if 
\begin{equation}\label{e101}
    \ell(p^{m+a}-k)\dfrac{\text{gcd}(d,24)^2}{\text{gcd}(d,24\ell)^2}-\ell p^{m-a}\dfrac{\text{gcd}(d,24p^a)^2}{\text{gcd}(d,24\ell)^2}+k \geq 0.
\end{equation}
To check the positivity of \eqref{e101}, we have to find all the possible divisors of $576\ell$.

We define $$\mathcal{S}=\{ 2^{r_1}3^{r_2}tp^s : 0 \leq r_1 \leq 6, 0\leq r_2 \leq 2, t \mid \ell \hspace{1.5mm} \text{but} \hspace{1.5mm} p\nmid t, \hspace{1.5mm} \text{and}\hspace{1.5mm} 0 \leq s \leq a \}$$  to be the set of all divisors $d$ of $576\ell$. We observe that the the values of $\dfrac{\text{gcd}(d,24)^2}{\text{gcd}(d,24\ell)^2}$ and $\dfrac{\text{gcd}(d,24p^a)^2}{\text{gcd}(d,24\ell)^2}$ are $\dfrac{1}{t^2p^{2s}}$ and $\dfrac{1}{t^2}$ respectively when divisors $d \in \mathcal{S}$. Now, substituting these values in the left side of \eqref{e101}, we see that it equals
\begin{equation}\label{e102}
     \dfrac{\ell}{t^2}\left( \dfrac{p^{m+a}-k}{p^{2s}} - p^{m-a} \right) + k.
 \end{equation}
When $s=a$, \eqref{e101} becomes $\left( k-\dfrac{k\ell}{t^2p^{2a}} \right) \geq 0$, which is true as $p^{2a}\geq \ell$. For $0 \leq s < a$, the quantity in \eqref{e102} is greater than equal to zero as $m>a$ and $k \leq p^{m+a}(1-p^{2s-2a})$ (by assumption). Hence the quantity in \eqref{e102} are greater than or equal to $0$ when $p^{2a} \geq \ell$ and $k \leq p^{m+a}(1-p^{2s-2a})$.

Therefore, the orders of vanishing of $\mathcal{B}_{\ell, p, m}(z)$ at the cusp $\dfrac{c}{d}$ is nonnegative. So $\mathcal{B}_{\ell, p, m}(z)$ is holomorphic at every cusp $\dfrac{c}{d}$. We have also verified the Nebentypus character by Theorem \ref{thm_ono1}. Hence $\mathcal{B}_{\ell, p, m}(z)$ is a modular form of weight $\dfrac{p^m(p^a-1)}{2}$ on $\Gamma_0(576\cdot \ell)$ with Nebentypus character $\chi_1(\bullet)$.
\end{proof}
 
	\begin{proof}[Proof of Theorem \ref{thm1}]Suppose $m>1$ is a positive integer. From Lemma~\ref{lem1}, we have 
		$$	\mathcal{B}_{\ell, p, m}(z)=\dfrac{\eta(24\ell z)^k\eta(24z)^{p^{m+a}-k}}{\eta(24p^a z)^{p^m}} \in M_{\frac{p^m(p^a-1)}{2}}\left(\Gamma_0(576 \cdot \ell), \chi_1(\bullet)\right),$$  Also the Fourier coefficients of the eta-quotient $\mathcal{B}_{\ell, p, m}(z)$ are integers. So, by Theorem \ref{serre} and Lemma \ref{lem2}, we can find a constant $\alpha>0$ such that
		$$
		\# \left\{n\leq X:T_{\ell, k}(n)\not\equiv 0 \pmod{p^m} \right\}= \mathcal{O}\left(\dfrac{X}{\log^{\alpha}{}X}\right).
		$$
Hence  $$\lim\limits_{X\to +\infty}\dfrac{	\# \left\{n\leq X:T_{\ell, k}(n)\equiv 0 \pmod{p^m} \right\}}{X}=1.$$ This allows us to prove the required divisibility by $p^m$ for all $m>a$ which trivially gives divisibility by $p^m$ for all positive integer $m\leq a$. This completes the proof of Theorem \ref{thm1}. 
 \end{proof}

\section{Concluding Remarks}\label{sec:conc}

\begin{enumerate}
    \item There are several more isolated congruences that we have found, which we do not mention in this paper. For instance, for all $n\geq 0$ it is easy to prove, using the theory of modular forms, that
    \[
    T_4(3^4n+57)\equiv 0 \pmod{12}.
    \]
    It would be interesting to make some more systematic studies and prove more general Ramanujan-type that the $T_{\ell, k}(n)$ function satisfies.
    \item There seems to be an infinite family of congruences for $T_9(n)$, similar to a result of Baruah and Das \cite{BaruahDas}. We leave it as an open problem to find such families.
    \item In Theorem \ref{thm1}, the bound on $k$ is given by $k \leq p^{m+a} \left( 1 - p^{2s-2a} \right) + \epsilon$, where we believe that there exists a positive integer $\epsilon$, although its exact value is yet to be determined. We leave this as an open problem.
    \item We close the paper with the following conjecture. 
\begin{conjecture}\label{conjp}
    Let $p\geq 5$ be a prime with $\left( \dfrac{-2}{p} \right)_{L}=-1$, and let $t$ be a positive integer with $\gcd(t,6) = 1$ and $p\mid t$. Then for all $n\geq 0$ and $1\leq j \leq p-1$, we have 
    \begin{equation*}
        T_2\left( 9\cdot t^2n + \frac{9\cdot t^2j}{p} + \frac{57\cdot t^2-1}{8} \right) \equiv 0 \pmod{6}.
    \end{equation*}
\end{conjecture}
Recently, Paudel, Sellers, and Wang \cite{PaudelSellersWang} have announced a proof of a generalization of Conjecture \ref{conjp}.
\end{enumerate}

\section*{Acknowledgements}
The authors are grateful to the anonymous reviewer for a very thorough reading of the manuscript and several helpful comments which improved the manuscript. The second author was partially supported by a Start-Up Grant from Ahmedabad University (Reference No. URBSASI24A5). The third author was partially supported by an institutional fellowship for doctoral research from Tezpur University, Assam, India.

\end{document}